\documentclass[journal]{IEEEtran}

\usepackage [english]{babel}

\IEEEoverridecommandlockouts

\usepackage[pdftex]{graphicx}
\usepackage{epstopdf}
\usepackage{amsmath}
\usepackage{amssymb}
\usepackage{amsfonts}
\usepackage{balance}

\usepackage{cite}

\usepackage{amsthm}%
\usepackage{balance}

\newtheorem{assumption}{Assumption}
\newtheorem{theorem}{Theorem}
\usepackage{amssymb}%

\usepackage[dvipsnames]{xcolor}

\usepackage{array}
\usepackage{multirow}
\usepackage{color}
\newcommand{\JW}{\textcolor{black}}
\newcommand{\il}{\textcolor{black}}
\newcommand{\jw}{\textcolor{black}}
\newcommand{\jdw}{\textcolor{black}}
\newcommand{\JDW}{\textcolor{black}}
\newcommand{\JDWW}{\textcolor{black}}
\newcommand{\ill}{\textcolor{black}}
\newcommand{\ak}{\textcolor{black}}
\newcommand{\an}{\textcolor{black}}
\newcommand{\ank}{\textcolor{black}}
\newcommand{\JWW}{\textcolor{black}}
\newcommand{\WJ}{\textcolor{black}}
\newcommand{\ka}{\textcolor{black}}
\newcommand{\icl}{\textcolor{black}}
\newcommand{\fv}{\textcolor{black}}
\usepackage[prependcaption,colorinlistoftodos]{todonotes}

\begin{document}

\title{\JWW{\fv{Frequency and voltage regulation} in hybrid AC/DC networks}}

\author{Jeremy D.~Watson and Ioannis Lestas
\thanks{{Jeremy D.~Watson and Ioannis Lestas are with the Department of Engineering, University of
Cambridge, Trumpington Street, Cambridge, CB2 1PZ, United Kingdom; emails: jdw69@cam.ac.uk, icl20@cam.ac.uk.}}
}
\maketitle

\begin{abstract}

Hybrid AC/DC networks are a key technology \JW{for future electrical power systems, due to the increasing number of converter-based loads and distributed energy resources.} In this paper, we consider the design of \il{control schemes} for hybrid AC/DC networks, focusing especially on the control of the interlinking converters (ILC(s)).
\icl{We present} two \il{control schemes}: firstly for decentralized primary control, and secondly, a \ank{distributed} secondary controller. In the primary case, the stability of the controlled system is proven in a general hybrid AC/DC network which may include asynchronous AC subsystems.
\an{Furt\JWW{h}ermore, it is demonstrated} that power-sharing across the AC/DC network is significantly improved compared to {\il{previously}} proposed dual droop control.
\WJ{The proposed scheme for secondary} \il{control} guarantees the convergence of the AC system frequencies \il{and the average DC voltage of each DC subsystem to their nominal values respectively}. \WJ{An optimal power allocation is also achieved at steady-state}.
 The \ak{applicability and \JDWW{effectiveness} of the proposed} algorithms are verified by simulation on a test hybrid AC/DC network in MATLAB / Simscape Power Systems.

\end{abstract}


\IEEEpeerreviewmaketitle

\section{Introduction}

\emph{\ak{Motivation:}} Hybrid AC/DC networks are a key technology \JW{for future electrical power systems.} One major reason is the increasing number of \JW{converter-interfaced energy sources and loads}. Direct current grids have several advantages \cite{jovcic2011} over traditional AC systems: lower power losses, largely due to the absence of reactive power; higher power transfer capability; and DC grids can also facilitate the connection of asynchronous AC grids. However, AC technology is well established and is more suitable for some applications. Therefore, it is likely that AC and DC technology will be combined via interlinking voltage source converters (VSC) to form a hybrid AC/DC network \cite{wang2017}.

Hybrid AC/DC networks present new challenges in terms of frequency and voltage control \cite{malik2017}.
In particular, \ak{an open problem} \cite{unamuno20152} is the control of the interlinking converter (ILC), \WJ{where the aim is to guarantee stability while ensuring that the frequency and voltages are appropriately regulated.}
\ak{This is a challenging problem since the ILC operation simultaneously affects the AC frequency and the DC voltage.  Moreover, a prescribed allocation is desired in many cases, such that economic optimality is achieved among generating units.} \fv{Furthermore,} \ak{distributed techniques for generation control are desirable due to the increasing penetration of renewable sources of generation which significantly increases the number of active elements in power grids, making \JWW{traditional} centralized approaches impractical and costly.}


\emph{\ak{Related work: }}\JDW{Numerous controllers for either AC or DC networks or microgrids alone have been proposed recently, e.g. from simple droop-based strategies to sliding mode control for DC networks in \cite{cucuzella2018}, distributed consensus for DC networks in \cite{morstyn2016,nasirian2015}, and model predictive control in \cite{papangelis2017}. For AC microgrids the literature is even more extensive, as surveyed in \cite{rocabert2012} for example. However, the control of a hybrid network presents new challenges as control actions on either the AC or DC sides affect the entire network. We therefore focus on the regulation of both AC and DC subsystems in a hybrid microgrid.}

Decentralized primary control strategies have significant advantages over centralized approaches, including: simplicity, plug-and-play capability, and immunity to communication failure \cite{unamuno20152}. The most common example, droop control, is cheap, intuitive, and easy to implement. AC frequency droop is ubiquitous, while DC voltage droop controllers are prevalent in \il{the} literature as well as in the few industrial multiterminal DC (MTDC) projects \cite{spallarossa2014}. The DC bus voltage dynamics are comparable to traditional AC frequency / real power control, in that an excess of power results in the MTDC voltages rising, while a shortage causes the MTDC voltages to decrease.

Therefore, an obvious way to control the interlinking converter is a dual droop scheme combining the two characteristics \il{in one} controller. The DC voltage droop stabilizes the DC grid and the MTDC systems participate in the frequency regulation of the connected AC grids \cite{chaudhuri2013}. However, the two \ak{droop schemes} interact with each other in a way that degrades their performance, as analyzed by \cite{akkari2015}, where a modified droop coefficient was calculated to prioritize the AC-side performance. \JWW{A strategy using the ILC power to improve performance is presented in \cite{malik2018}, although the coupling between the AC and DC grids still introduces some inaccuracy.} \JW{Another approach is seen in \cite{papangelis2017epsr} which uses receding horizon control to \il{provide a compromise} between AC frequency and DC voltage deviations.}

\jw{A new idea for the control of the ILCs was proposed in \cite{wang2017}, \cite{loh2013}, and \cite{luo2016}. The per-unit values of the AC frequency and DC voltage are synchronized by controlling the power transfers with a proportional-integral controller. This allows the ILC to perform network emulation by relating the AC frequencies and DC voltages to each other; and power sharing is then set by the droop coefficients on both the AC and DC sides.}

Secondary controllers with communication aim to restore the frequencies and voltages to their nominal value and to share power equitably between sources in the AC and DC subsystems \cite{malik2017}. In \cite{dai2011}, a distributed controller for sharing frequency reserves of asynchronous AC systems via HVDC \il{was} designed. However the DC voltage dynamics were not modelled. \JW{The authors in \cite{andreasson2017} designed distributed controllers for distributed frequency control of asynchronous AC systems connected through a MTDC grid.}

\emph{\ak{Contribution: }}
\JW{In this paper, we present new voltage source converter control schemes for the interlinking converters which are inspired by the controller proposed in \cite{jouini2016}.} The proposed \ak{controllers use} the energy stored in the DC capacitance as the "inertia" for the emulation of synchronous machines. \JW{We show that this idea can successfully be applied to the interlinking converter in general hybrid AC/DC networks for both primary and secondary control, \ak{and that the proposed schemes have} advantages over previously proposed ILC controllers.}
\ak{In particular, our decentralized control approach, by making use of the energy stored in DC-side capacitance,
achieves faster primary frequency} regulation  \jw{compared to schemes that directly control the power transfer}.
\ak{Moreover, we provide stability guarantees for primary frequency regulation and sufficient conditions for optimal steady-state power allocation.}

 \WJ{We also propose a \fv{scheme} for secondary frequency and voltage \fv{control} which}
\icl{ regulates the frequency and the weighted average voltages of DC subgrids to  prescribed nominal values at steady state, thus providing improved frequency and voltage regulation, and also leads to an optimal power sharing.} Moreover, we show that virtual capacitance in the controller can be used \ak{to   further improve} performance.

\ak{For clarity, we summarize the main contributions of the paper below:
\begin{enumerate}
\item We \WJ{propose a} decentralized VSC controller \WJ{inspired by} \cite{jouini2016} for general hybrid AC/DC networks. For this setting, we provide stability guarantees and sufficient conditions for an optimal steady state power allocation.
    \item We propose a novel approach for the control of interlinking converters and generation sources \ka{for secondary frequency and voltage regulation} \il{which} guarantees the convergence \il{of both the AC system frequencies and the weighted average DC voltage of each DC subsystem} to their nominal values. \WJ{A prescribed power sharing is also achieved. We also show that virtual capacitance in the controller can be used to adjust the \fv{voltage deviations} and improve performance.}
\end{enumerate}
}

\emph{\ak{Paper structure: }} The paper is organized as follows.
Section \ref{formulation} presents the network \JW{model} and formulates the control problem. Our main results are given in section \ref{results}, including the proposed primary controller in section \ref{pri} and the \ank{distributed} secondary controller in section \ref{sec}. The performance of the two controllers is illustrated via \JDW{case studies} in section \ref{casestudy} and compared to traditional controllers. Finally, conclusions are presented in section \ref{conclusion}. \ak{The proofs of all the results presented can be found in the appendix.}

\section{Problem formulation}\label{formulation}

\subsection{Network model}

\begin{figure}[!ht]
    \centering
    \includegraphics[width=0.45\textwidth]{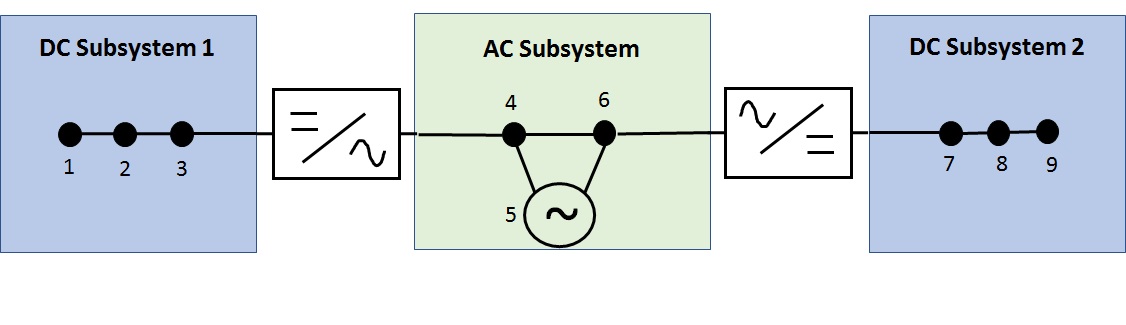}
    \caption{\il{AC/DC network diagram.}}
    \label{fig:network}
\end{figure}

\JW{We consider a general hybrid AC/DC network \il{with the set of buses denoted by} $N = (1,2, ..., |N|)$ and \il{the set of transmission lines by} $E = (1,2, ..., |E|)$. The network is composed of multiple AC and DC subsystems. \il{We denote the set of subsystems by $K = (1,2, ..., |K|)$ and we also have} $N = (\cup_i N_i^{dc}) \cup (\cup_j N_j^{ac})$, where $N_i^{ac}$ and $N_j^{dc}$ denote the collection of buses belonging to the AC subsystem $i \in K$ and DC subsystem $j \in K$ respectively. \il{Each subsystem is assumed to be connected and it is connected to the rest of the network only via \ak{interlinking converters}\footnote{\il{Note that this is without loss of generality since the connection of two collections of AC buses (or DC buses respectively) may simply be considered as one larger subsystem.}}, as} illustrated in Figure \ref{fig:network}. Each AC subsystem $i$ \il{is} described by the connected graph $(N_i^{ac}, E_i^{ac})$ with arbitrary direction, and each DC subsystem $j$ by the connected graph $(N_j^{dc}, E_j^{dc})$ with arbitrary direction\footnote{\jw{The results presented in the paper are unaffected by the choice of direction}.}. For each bus $j \in N$ we use \ak{$i : i\to j$ and $k : j\to k$} to denote the predecessors and successors of bus $j$ respectively. For convenience we also define the set of all AC buses $N_{ac} = \cup_{i} N_i^{ac}$ and all DC buses $N_{dc} = \cup_{j} N_j^{\WJ{d}c}$, such that $N_{ac} \cup N_{dc} = N$; likewise we define the set of all AC edges $E_{ac} = \cup_{i} E_i^{ac}$ and all DC edges $E_{dc} = \cup_{j} E_j^{dc}$\ak{.
Connections} between AC and DC buses are facilitated by the interlinking converters, the set of which is \il{denoted} by $X = (1,2, ..., |X|)$. \il{The ILC buses are denoted by $X_{\ka{j}}^{ac}\in N_{ac}$ and $X_{\ka{j}}^{dc}\in N_{dc}$ for the AC and DC buses, respectively,} to which the ILC $\ka{j}$ is connected. \il{The set of AC buses \il{to which a converter is connected} \il{is} denoted by \il{$X_{ac}= (X_1^{ac},\hdots,X_{|X|}^{ac}) \subset N_{ac}$}.}} \jw{Similarly, the set of DC buses to which a converter is connected is denoted by $X_{dc}= (X_1^{dc},\hdots,X_{|X|}^{dc}) \subset N_{dc}$.}
\begin{table}[!ht]\centering
\caption{\il{Notation in system model}}
\renewcommand{\arraystretch}{1.25}
\begin{tabular}{ c l }
 \hline
 $\omega_j$ & AC frequency deviation at bus $j$\\
 $\eta_{ij}$ & voltage angle \il{between} two AC buses $i$ and $j$ \\
 $M_j$ & inertia at bus $j$ \\
 $p^G_j$ & \jw{generated power at bus $j$} \\
 $p^L_j$ & \il{load power} at bus $j$ \\
 $p^X_j$ & interlinking converter \il{power transfer} at bus $j$ \\
 $p_{ij}$ & \il{power transfer} between buses $i$ and $j$ \\
 $B_{ij}$ & \il{transmission line susceptance} for $(i, j) \in E_{ac}$ \\
 $D_j$ & damping coefficient at bus $j \in N_{ac} \setminus X_{ac}$ \\
 $C_{j}$ & capacitance at bus $j \in N_{dc}$ \\
 $V_{j}$ & DC voltage \jw{deviation} at bus $j$ \\
 $G_{ij}$ & line conductance for $(i, j) \in E_{dc}$\\
 \hline
\end{tabular}
\label{tab:not}
\end{table}
\begin{assumption}\label{as:1}
We make the following assumptions regarding the network:
\begin{itemize}
    \item 1a: Voltage magnitudes are 1 p.u. for all buses $j \in N_{ac}$.
    \item 1b: Lines $(i,j) \in E_{ac}$ are lossless and are characterized by their constant reactances $X_{ij} > 0$.
    \item 1c: Reactive power does not \ak{affect} either bus voltage angles or the frequency, and is thus ignored.
    \item 1d: The AC system(s) are three-phase balanced.
    \item 1e: Bus voltages are close to 1 p.u. for all $j \in N_{dc}$, \jdw{such that currents and powers are approximately equivalent in a per-unit system.}
    \item 1f: Lines $(i,j) \in E_{dc}$ are characterized by their conductance \jdw{$G_{ij} = \frac{1}{R_{ij}}$, where $R_{ij}$ is the line resistance.} \JDW{The line losses are small and may be neglected.}
\end{itemize}
\end{assumption}

\begin{figure}[!ht]
    \centering
    \includegraphics[width=0.4\textwidth]{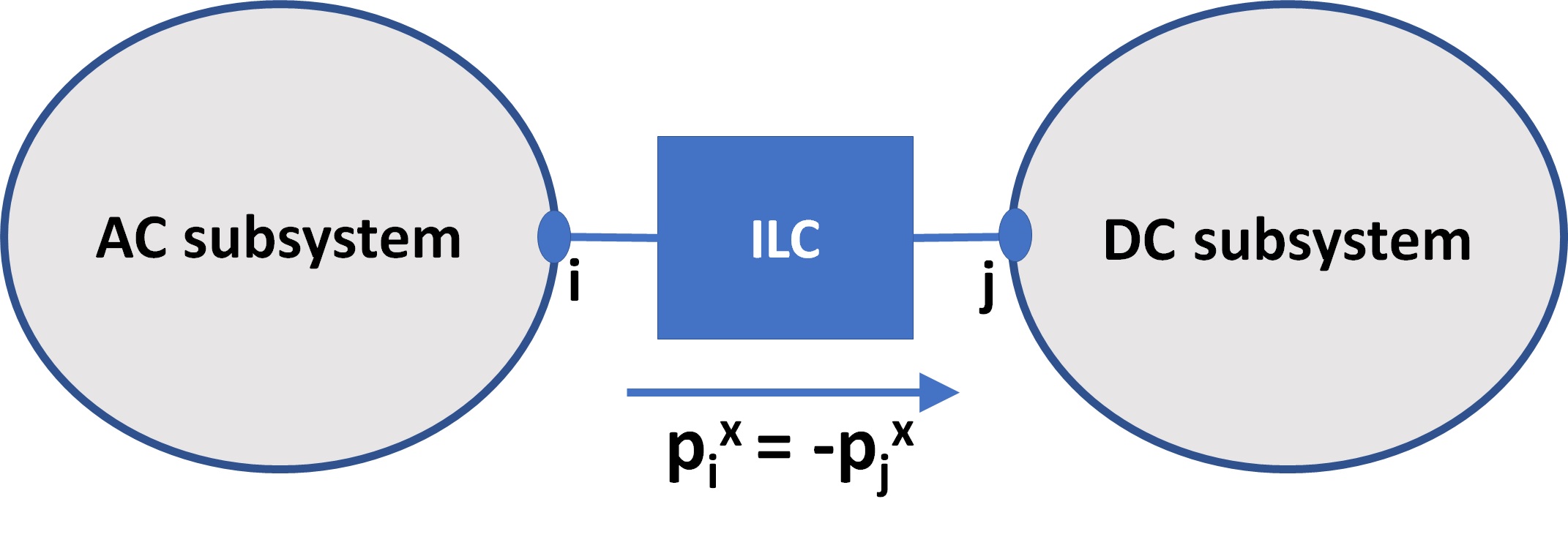}
    \caption{ILC connection diagram}
    \label{fig:ilc}
\end{figure}

\emph{Remark 1:} Assumption \ref{as:1} may be explained as follows:
\begin{itemize}
    \item 1a-d: These are well-known assumptions for AC transmission systems which are used in much of the literature. These assumptions allow us to model the active power transfer through a line $(i,j)$ as $p_{ij} = B_{ij}$ sin $\eta_{ij}$ where $B_{ij} = 3X_{ij}^{-1} > 0$.
    \item 1e-f: These are typical assumptions \jw{in DC networks} \JDW{\cite{andreasson2017}}.
\end{itemize}
We also consider the modelling of the interlinking converters, as illustrated in Fig. \ref{fig:ilc}. The AC-side \il{bus} of the interlinking converter has no inertia of its own, however power imbalances \il{in the AC subsystem lead to a power transfer through the converter. Note also that this power transfer will affect the DC-side voltage of the converter.} \jw{For each bus $j \in X_{ac} \cup X_{dc}$, we define the power transfer $p_j^X$ as the power leaving the bus through the ILC. Hence for a converter bus $j \in X_{ac}$, the power transfer $p_j^X$ is the AC-to-DC transfer, whereas for a converter bus $j \in X_{dc}$, the power transfer $p_j^X$ is the DC-to-AC transfer. We assume that such power transfers are instantaneous and lossless, hence for an ILC $x$ we have \jdw{$p_{X_x^{dc}}^X=-p_{X_x^{ac}}^X$.}}
 Given these assumptions and definitions, the hybrid AC/DC network dynamics are:
\begin{subequations}\label{eq:system}
\begin{align}
    \dot{\eta}_{ij} &= \omega_i - \omega_j, (i, j) \in E_{ac}\\
    M_j\dot{\omega}_j &= p_j^G - p_j^L + \ak{\sum_{i : i\rightarrow j}}p_{ij} - \ak{\sum_{k : j\rightarrow k}}p_{jk} -D_j\omega_j,\nonumber\\ &\ \ \ \ j \in N_{ac} \setminus X_{ac}\label{eq:swing}\\
    0 &= p_j^G - p_j^L + \ak{\sum_{i : i\rightarrow j}}p_{ij} - \ak{\sum_{k : j\rightarrow k}}p_{jk}- p_j^X, \an{j \in X_{ac}}\label{ref:eq:convbus}\\
    C_j\dot{V}_j &= p_j^G - p_j^L + \ak{\sum_{i : i\rightarrow j}}p_{ij} - \ak{\sum_{k : j\rightarrow k}}p_{jk} - \jw{p_j^X}, \an{j \in N_{dc}}\label{eq:dcv}\\
    p_{ij} &=
    \begin{cases}
    B_{ij}\sin{\eta_{ij}} & j \in N_{ac}\\
    G_{ij}(V_i-V_j) & j \in N_{dc}
    \end{cases}, (i, j) \in E \label{eq:pf}
\end{align}
\end{subequations}
We now write the system dynamics in matrix form. \JW{The vector of angle differences is $\eta = [\eta_{ij}]_{(i,j)\in E_{ac}}$, the vector of \il{AC frequency deviations from its nominal value (50 or 60~Hz)} is denoted by $\omega = [\omega_j]_{j \in N_{ac}}$, while the vector of DC \il{voltage deviations from their nominal value} is denoted by $V = [V_j]_{j\in N_{dc}}$. $M$ is the diagonal matrix of the generator inertias $M_j, j \in N_{ac} \setminus X_{ac}$, while the damping coefficients are $D = diag([D_j]_{j \in N_{ac} \setminus X_{ac}}$. The frequencies at the corresponding AC buses are denoted by $\omega^G = \ak{[\omega_j]_{j \in N_{ac} \setminus X_{ac}}}$, and the vector of frequencies at the converter buses is \il{$\omega^X = [\omega_j]_{j \in X_{ac}}$}. The diagonal matrix of the DC bus capacitances is $C = diag([C_j]_{j \in N_{dc}})$. The vector of generator powers is denoted by $p^G = [p_j^G]_{j \in N}$, the vector of load powers by $p^L = [p_j^L]_{j\in N}$. \jdw{We also use the notation $p_j^X = 0$ for buses without converters, i.e. $j \in N \setminus (X_{ac} \cup X_{dc})$ and denote the vector of converter powers by $p^X = [p_j^X]_{j \in N}$.} \JWW{Similarly, at the converter buses $j \in X_{ac}$ we use the notation $p_j^G = 0$.} The power transfer vector is defined by $p^F = [p_j^F]_{j \in N}$, where each $p_j^F = \sum_{j : i\rightarrow j}p_{ij} - \sum_{j : j\rightarrow k}p_{jk}$. The matrix $A$ is the incidence matrix of the graph $(N_{ac}, E_{ac})$.} The system equations are thus:
\begin{subequations}
\label{ref:sys}
\begin{align}
    \dot{\eta} &= A^T\omega= A^T \begin{bmatrix}
    \omega^G \\
    \omega^X
    \end{bmatrix}\\
    \begin{bmatrix}
    M \dot{\omega}^G\\
    0 \\
    C\dot{V}
    \end{bmatrix}
    &= p^G-p^L-p^X-p^F-
    \begin{bmatrix}
    D\omega^G \\
    0\\
    0
    \end{bmatrix}
\end{align}
\end{subequations}
\subsubsection{Equilibrium conditions}
An equilibrium of the system in (\ref{ref:sys}) is defined by the following conditions:
\begin{subequations}
\label{eqconditions}
\begin{align}
    0 &= A^T\omega= A^T \begin{bmatrix}
    \omega^G \\
    \omega^X
    \end{bmatrix}\\
    0 &= p^G-p^L-p^X-p^F-
    \begin{bmatrix}
    D\omega^G \\
    0\\
    0
    \end{bmatrix}
    \end{align}
\end{subequations}

\JW{
\ak{We assume that there exists\footnote{\icl{Existence of equilibria in AC systems is beyond the scope of this paper and have been considered in e.g. \cite{dorfler2013synchronization}.}} some equilibrium point of \eqref{ref:sys},} and denote such an equilibrium by $(\eta^\ast, \omega^{G\ast}, V^\ast)$.
Individual equilibrium values are also denoted by the superscript asterisk, e.g. $\eta^\ast_{ij}$, $\omega_i^{G\ast}$, and $V_j^\ast$.}\\
\begin{assumption}\label{as:2}
  \JW{$|\eta^\ast_{ij}| < \frac{\pi}{2}$ for all $(i, j) \in E_{ac}$}.
\end{assumption}

\il{This assumption is often referred as a security constraint and is common in the literature for power grid stability analysis.}
\subsection{Control objectives}

The control objectives are:
\begin{enumerate}
    \item Solutions must converge \il{to an equilibrium point}.
    \item \JDW{For primary control, AC frequencies and DC voltages should not deviate too far from their nominal values, i.e. $\lim_{t \to \infty} |\omega_j(t)| < e_\omega$ for all buses $j \in N_{ac}$ and $\lim_{t \to \infty} |V_j(t)| < e_V$ for all buses $j \in N_{dc}$ for some appropriate scalars $e_V$ and $e_\omega$. }
    \item AC frequencies \WJ{and \fv{a} weighted average of the DC voltages should} converge to their nominal values for secondary control.
    \item \JDW{Power sharing between all sources should be optimal.}
\end{enumerate}

The last objective may be stated more formally by considering the minimization of a quadratic cost function \cite{monshizadeh2017}:
\begin{subequations}\label{eq:optim}
\begin{align}
    \underset{p^G}{\text{min}} \mbox{   } C_G &= \frac{1}{2}(p^G)^TQp^G\label{eq:cost}\\
    \text{subject to:} \mbox{   }
    \mathbf{1}^Tp^G &= \mathbf{1}^Tp^L + \mathbf{1}^T\begin{bmatrix}
    D\omega^G \\
    0\\
    0
    \end{bmatrix}\label{eq:constraint} \\
    \an{p^G_j = 0, j \in X_{ac}} \label{gen:constraint}
\end{align}
\end{subequations}
\ak{where} $Q$ is a \ak{positive definite} diagonal matrix containing the cost coefficients for each energy resource, \il{and $\mathbf{1}$ is the vector of ones \ka{with appropriate dimension}.}
\an{Note that constrain\JWW{t} \eqref{eq:constraint} is a requirement \ka{for} power balance at equilibrium, i.e. that the total generation and demand are equal, while \eqref{gen:constraint} suggests that the generation at converter AC buses is zero, \ka{which holds by definition  (note that $p^G_j$ appears in \eqref{ref:eq:convbus} only for convenience in presentation).} }
\an{To proceed further, we define the diagonal matrix $\tilde{Q}$ such that $\tilde{Q}_{ii} = Q_{ii}^{-1}, i \in N / X_{ac}$ and $\tilde{Q}_{ii} = 0, i \in X_{ac}$. With slight abuse of terminology, we shall refer to $\tilde{Q}$ as the inverse cost matrix.}
%
 Using the standard method of Lagrange multipliers as in \cite{monshizadeh2017}, and defining the vector $p^u = [D(\omega^G)^T \mbox{ }0^T\mbox{ } 0^T]^T$ for convenience, the solution $p^{G\ast}$ to the optimization problem is:
\begin{equation}\label{eq:optimsoln}
    p^{G\ast} = \frac{\JWW{\tilde{Q}}\mathbf{1}\mathbf{1}^T}{\mathbf{1}^T\JWW{\tilde{Q}}\mathbf{1}}(p^L+p^u)
\end{equation}

\section{Main results} \label{results}

\subsection{Decentralized primary control} \label{pri}

We assume power-frequency droop control for the AC generators and power-voltage droop for the DC energy resources:
\begin{equation} \label{eq:droop}
    p^G = -\an{\tilde{Q}}
    \begin{bmatrix}
    \omega \\
    mV
    \end{bmatrix} + \jw{p^G_{nom}}
\end{equation}
\ak{where} $\an{\tilde{Q}} \geq 0$ is the inverse cost matrix \jw{of droop coefficients}, and $m > 0$ is a constant, and as stated previously, $\omega$ and $V$ are the column vectors of the AC \il{frequency and DC voltage deviations,} respectively. \jdw{The nominal power generation $p^G_{nom}$ is the power allocation that satisfies \eqref{eq:optim} when the frequencies are at their nominal values.} \JDW{The second control objective (limitation of frequencies and voltages deviations) may be satisfied by choosing suitably large droop coefficients in $\an{\tilde{Q}}$.} In order to simplify the presentation here we use \il{proportional droop control schemes}, nevertheless this condition could be relaxed to local input strict passivity of the dynamics of each AC generator from input $-\omega_j$ to output $p_j^G$ and each DC generator from input $-V_j$ to output $p_j^G$ around their respective equilibrium values $\omega_j^\ast$ and $V_j^\ast$, similarly to the \ka{analysis} in \cite{kasis2017, Power3b, kasis2017stability}, .

\JDW{We also introduce} a VSC controller based on \cite{jouini2016}. Let the voltage \il{angle $\theta_i$} at the AC-side output of an ILC $x$ be:
\begin{subequations} \label{eq:ilcpri}
\begin{align}
    \theta_i &= \int{mV_j}, \  \text{\il{ i.e. the frequency is given by,}}\\
    \omega_i &= \dot{\theta}_i = mV_j \label{eq:ilc}
\end{align}
\end{subequations}
\ak{where} $i \in X_x^{ac}$ and $j \in X_x^{dc}$. \JW{This relates the AC frequency deviation proportionally to the DC voltage deviation by a chosen constant $m > 0$.}
 We assume that the converter is lossless and that the internal dynamics are sufficiently fast compared to the network dynamics. In \cite{jouini2016} the suggestion is to set $m = \frac{\omega^{nom}}{V_{dc}^{nom}}$ where $\omega^{nom}$ and $V_{dc}^{nom}$ are the nominal values of the AC grid frequency and the DC grid voltage. \JW{Since in this paper $\omega_i$ and $V_j$ are deviations from a nominal value, other values of $m$ are also possible. Large values of $m$ result in smaller DC voltage deviations and larger AC frequency deviations, and in general $m = \frac{\omega^{nom}}{V_{dc}^{nom}}$ may be too large for this \il{scheme} as frequency deviations are generally less acceptable than voltage deviations.}

Instead of directly controlling the power transfer, (\ref{eq:ilc}) \il{relates the frequency and voltage within the AC and DC sides, respectively}. Not only does this improve the accuracy of the power-sharing in the primary time-frame, but \ka{also provides} fast response to AC disturbances via capacitive inertia as discussed in \cite{jouini2016} and \cite{monshizadeh2017}. \jw{\ak{In} this paper, we take this concept further and use the capacitive inertia of the entire DC subsystem for frequency support, and also use the inertia of the AC system to regulate the DC voltage when appropriate.}

\begin{theorem}[Stability]\label{conv_thm}
  \JW{Consider a dynamical system described \il{by equations \ill{\eqref{eq:system}, \eqref{ref:sys}} with the control scheme} in \ka{(\ref{eq:droop}), (\ref{eq:ilc}),} \ill{and an equilibrium point for which Assumption~\ref{as:2} holds. Then \ak{there exists an open neighbourhood of this equilibrium} point such that all solutions of
\ak{\eqref{eq:system},\eqref{ref:sys},\eqref{eq:droop},\eqref{eq:ilc} starting in this region} converge to the set of equilibrium points as defined in (\ref{eqconditions}).}}
\end{theorem}

\ak{Theorem \ref{conv_thm} demonstrates the local convergence of solutions to \eqref{eq:system},\eqref{ref:sys},\eqref{eq:droop},\eqref{eq:ilc} to the set of its equilibria. Note that the result is local due to the sinusoidal power transfers in \eqref{eq:pf} and that it becomes global if those are linearized.}

\ak{The following theorem demonstrates that when line resistances become arbitrarily small, then the equilibria of the \ka{considered} system tend towards the global minimum of \eqref{eq:optim}.
}

\begin{theorem}[Power sharing]\label{optim_thm}
\jw{As the line resistances become arbitrarily small, the power sharing of \ka{the system \eqref{eq:system}, \eqref{ref:sys} with} the control scheme (\ref{eq:droop}), (\ref{eq:ilc}) becomes arbitrarily close to the solution of the optimization problem \eqref{eq:optim}.}
\end{theorem}

\emph{Remark 2:} In a practical network there will always be some small line resistances which affect power sharing. Nevertheless, if these resistances are small, \WJ{the} voltage \WJ{deviations} will also be small and the power sharing will be close to optimal.

\emph{Remark 3 (Power sharing in a dual-droop ILC controller scheme): } \JW{Consider the \il{dual-droop scheme \eqref{eq:dualdroop} often used in the literature} for primary control of the ILC, }
\begin{equation}\label{eq:dualdroop}
    p_i^X = K_i^{\omega}\omega_i - K_j^V{V_j}
\end{equation}
\JW{\ak{where} $K_i^{\omega}$ and $K_j^V$ are the respective droop coefficients, and the power transfer is directly controlled\footnote{In practice, the ILC controls the power transfer by varying its output voltage angle until \eqref{eq:dualdroop} is satisfied.}. It is clear that \eqref{eq:dualdroop} is unable to guarantee correct power-sharing for a disturbance at any arbitrary bus under the same assumptions. For droop-controlled sources to contribute power in proportion to their droop coefficients, a system-wide synchronizing variable is required. The proposed controller \eqref{eq:ilc} achieves this by relating the AC frequency to the DC voltages. By contrast, the dual droop controller \eqref{eq:dualdroop} \jdw{does not provide such a relation.}}

\subsection{\ank{Distributed} secondary control} \label{sec}

In this section we propose a \ank{distributed} controller inspired by \cite{monshizadeh2017} and \cite{andreasson2016}.
The concept of \textit{network emulation} can be carried further with the aid of distributed communication. Let the average DC voltage deviation of each subsystem $k$ be represented by the capacitance-weighted average $\bar{V}_k$:
\begin{equation}\label{eq:barv}
    \bar{V}_k = \sum_{j\in N_k^{dc}}C_jV_j
\end{equation}

As in the second distributed MTDC controller proposed in \cite{andreasson2016}, the DC voltages within each subsystem are either \il{communicated within the network so as to obtain $\bar V_k$ (for small subsystems) or} \il{$\bar V_k$ is obtained via an appropriately fast distributed approach, such that its dynamics can be decoupled from the stability analysis in this paper.}
From (\ref{eq:barv}) we have:
\begin{equation}\label{barv1}
\begin{aligned}
    \dot{\bar{V}}_k &= \sum_{j\in N_k^{dc}}C_j\dot{V}_j = \jw{\sum_{j\in N_k^{dc}}(p_j^G-p_j^L-p_j^X+p_j^F)}
\end{aligned}
\end{equation}
The DC branch-flows $p_j^F$ cancel out within the subsystem, and we therefore have the following expression which resembles the swing equation:
\begin{equation}\label{barv2}
    \dot{\bar{V}}_k = \sum_{j\in N_k^{dc}}(p_j^G-p_j^L-\jw{p_j^X})
\end{equation}

We now introduce the concept of \textit{virtual frequency deviation} $\hat{\omega}$ which is defined for the entire network as follows:
\begin{equation}\label{eq:virtualw}
    \hat{\omega}_j =
    \begin{cases}
    \omega_j &\mbox{if } j \in N_{ac} \\
    m\bar{V}_k & \mbox{if } j \in N_k^{dc}
    \end{cases}
\end{equation}
where $k$ is the DC subsystem to which all nodes in the associated set $N_k^{dc}$ belong. We will denote the vector of average DC subsystem voltages by $\bar{V}$. The converter which \il{interlinks} AC bus $i$ and the DC subsystem $k$ is governed by
\begin{equation}\label{eq:secILC}
    \omega_{i}  = m\bar{V}_k
\end{equation}
\ak{where} $m > 0$ is a positive coefficient. A common approach \il{to achieve \ka{optimal} power sharing in secondary control} is to introduce a synchronizing communicating variable $\xi$, e.g. \cite{monshizadeh2017}, and update these values \JW{via distributed averaging through an undirected connected communication graph. In particular, we denote this graph by \an{$(N,\tilde{E})$, where $\tilde{E}$ denotes the set of communication links, \icl{and also} denote by $\mathcal{L}$ the Laplacian of $(N,\tilde{E})$, defined as
\begin{equation}\label{dfn_laplacian}
\mathcal{L}_{ij} = \begin{cases}
\icl{\text{deg}(i)}, \mbox{ if } i=j, \\
-1, \mbox{ if } (i,j) \in \tilde{E} \\
0, \mbox{ otherwise }
\end{cases}
\end{equation}
where \icl{$\text{deg}(i)$} denotes the degree of node $i$.}} Then the \ank{distributed} controllers for the hybrid AC/DC system are:
\begin{subequations}\label{eq:sec}
\begin{align}
     \ak{T_\xi}\dot{\xi} &= -\mathcal{L}\xi -
     \an{\tilde{Q}}\ak{\hat{\omega}} \label{eq:firstpg}\\
     p^G &= \an{\tilde{Q}}
     \xi\label{eq:secpg}
\end{align}
\end{subequations}
\an{where $T_\xi$ denotes the diagonal matrix with positive time constants}, \JW{$\xi$ is the column vector of the synchronizing variables $\xi_j$},
 and $\an{\tilde{Q}}$ is the inverse cost coefficient matrix as before.

 The DC bus voltages are weighted by the associated capacitances in order to capture the dynamics of the physical energy of the subsystem \il{as follows from \eqref{barv1}}. \ak{However, buses with low capacitance could} have voltages far from the nominal while still satisfying $\bar{V}_k = 0$. \ak{Nevertheless}, as the DC bus voltages must still satisfy the power-flow equations, the voltages of buses with low capacitance will still be close to nominal. Furthermore, virtual capacitance $C_j^V$ may easily be added at any DC source bus $j$ via a derivative term in the DC source dynamics, e.g.:
\begin{equation*}
     p_j^G = \an{\tilde{Q}}_{jj}\xi_j - C_j^V\dot{V_j}
\end{equation*}
\ak{noting that the addition of the derivative term will not affect the steady-state value of $p^G_j$, \ka{hence} allowing it to retain its optimality properties.}

\ak{Our first result, proven in the appendix, demonstrates that the introduction of the controller \eqref{eq:virtualw},\eqref{eq:secILC},\eqref{eq:sec} ensures that the equilibria of the system \eqref{ref:sys},\eqref{eq:virtualw},\eqref{eq:secILC},\eqref{eq:sec} coincide with the global minimum of the optimization problem \eqref{eq:optim}.}

\begin{theorem}[Power sharing] \label{optm_thm}
\il{An equilibrium of the system \ka{\eqref{eq:system}},\eqref{ref:sys} with the control scheme \eqref{eq:virtualw},\eqref{eq:secILC},\eqref{eq:sec} solves} the optimization problem \eqref{eq:optim}.
\end{theorem}

\ak{The following theorem, proven in the appendix, demonstrates the local convergence of solutions \an{of the dynamical system \eqref{eq:system},\eqref{ref:sys}, when the controller
\eqref{eq:virtualw},\eqref{eq:secILC},\eqref{eq:sec} is applied,} to the global minimum of the optimization problem \eqref{eq:optim}. Furthermore, it guarantees \an{that frequency returns to its nominal value at equilibrium}, i.e. that $\omega^\ast = 0$, and that the average voltage deviation in every DC subsystem is zero, i.e. that $\bar{V}^\ast = 0$.}

\begin{theorem}[Convergence to optimality]\label{conv_opt_thm}
 \JW{Consider the dynamical system described in \ill{\eqref{eq:system},\eqref{ref:sys}} with the \il{control scheme \eqref{eq:virtualw},\eqref{eq:secILC},\eqref{eq:sec}} \ill{and an equilibrium point for which \ak{Assumption} \ref{as:2} is satisfied. \ak{Then, there exists an open neighbourhood of this equilibrium point such that all
solutions of \eqref{eq:system},\eqref{ref:sys},\eqref{eq:virtualw},\eqref{eq:secILC},\eqref{eq:sec} starting in this region
 converge to a set of equilibria that solve the optimization problem \eqref{eq:optim}, with $\omega^\ast = 0$ and $\bar{V}^\ast = 0$.}}}
 \end{theorem}

Theorem \ref{conv_opt_thm} demonstrates that all solutions of the considered system locally converge to an optimal solution of \eqref{eq:optim}. 

\emph{Remark 4:} \icl{Our proposed controller is distributed in the sense that
its implementation in a DC subgrid makes use of voltage measurements only within that subgrid,  
and is also fully distributed in the AC subgrids of the network. It should be noted that relaxing  \eqref{eq:barv} to a fully distributed controller that makes use of only local voltage measurements \fv{without additional information transfer}, while retaining the stability and optimality properties presented, is a highly non-trivial problem as this would distort the synchronization of the communicating variable $\xi$ needed for optimal power sharing.}

\section{Case studies}  \label{casestudy}

\begin{figure}[!ht]
    \centering
    \includegraphics[width=0.5\textwidth]{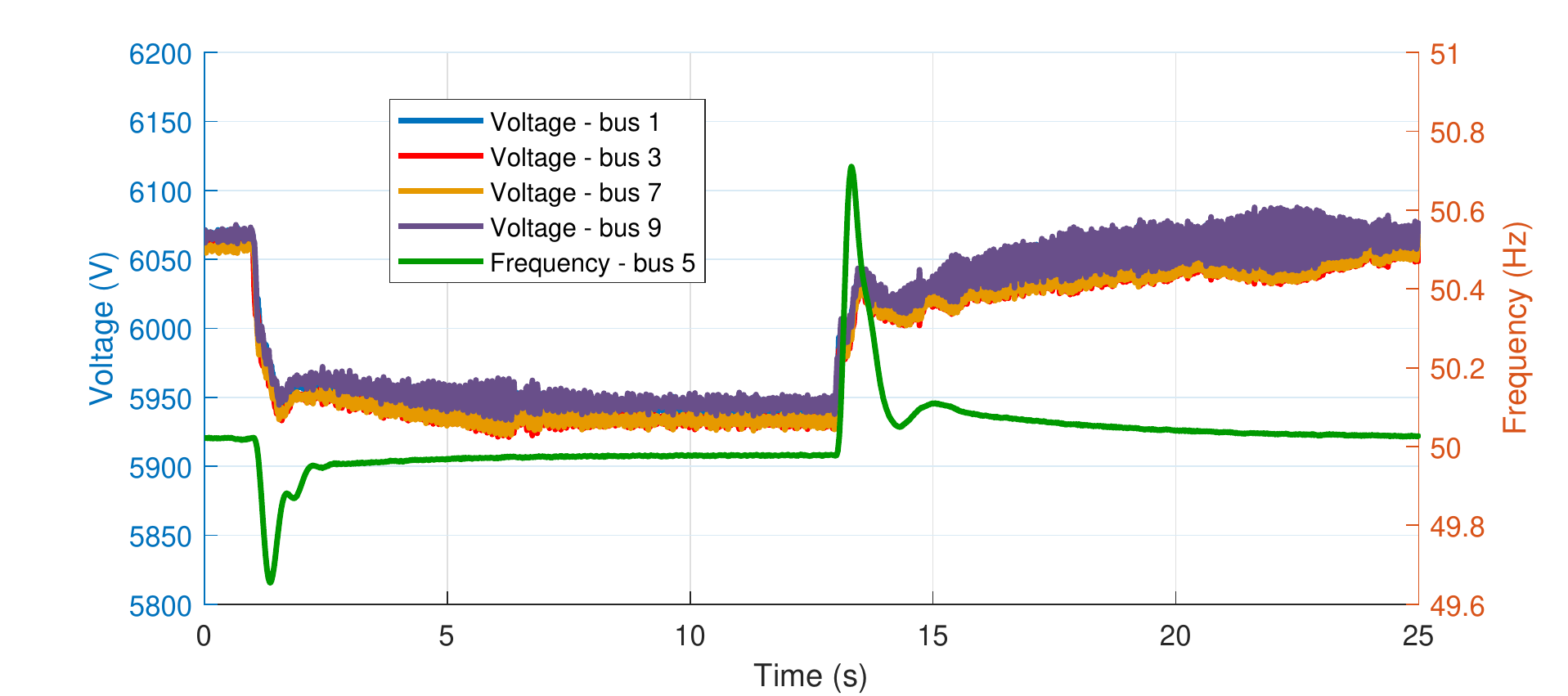}
    \caption{Frequency and voltage response with the decentralized primary controller \eqref{eq:droop},\eqref{eq:ilcpri}}
    \label{fig:b1fv}
\end{figure}

\begin{figure}[!ht]
    \centering
    \includegraphics[width=0.5\textwidth]{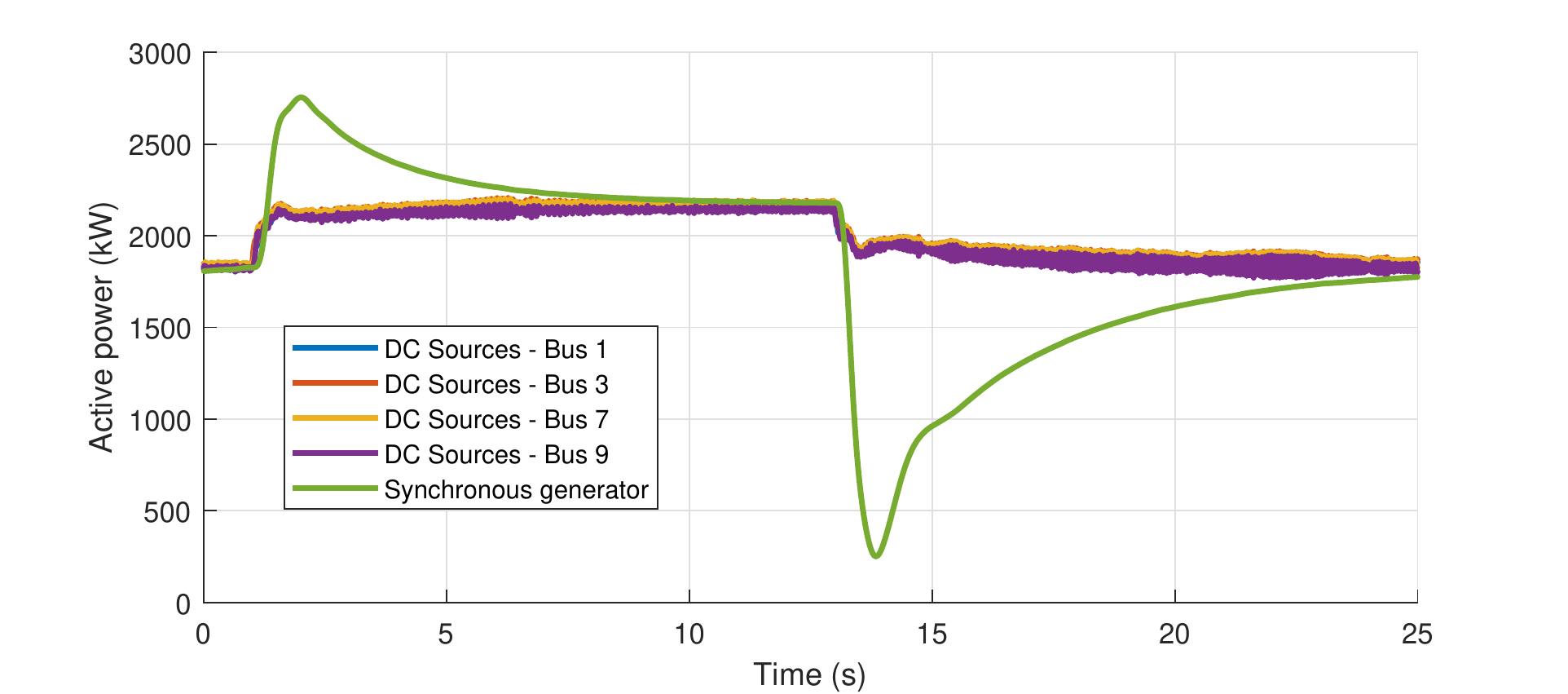}
    \caption{Power-sharing response with the decentralized primary controller \eqref{eq:droop},\eqref{eq:ilcpri}}
    \label{fig:b1ps}
\end{figure}

In order to demonstrate our results, \JDW{we study} the hybrid AC/DC system shown in Fig. \ref{fig:network} in MATLAB / Simscape Power Systems. The parameters of the network are given in Table \ref{tab:acdcparams}. \JDW{The synchronous machine is 4MVA, 13.8kV and is connected to the AC subsystem at bus 5 via a 4 MVA transformer, and there are eight distributed DC sources across the two subsystems at buses 1,3,7 and 9. All droop gains are chosen such that \JWW{the optimal power allocation is equal contributions from all sources.} 
The simulation model is more detailed and realistic than our analytical model, and it includes the inverter dynamics \ak{with} switching (two-level pulse width modulation), line dynamics, detailed generator, turbine-governor, and exciter dynamics, and realistic communication delays. In the simulation, we obtain $\bar{V}_k$ in \eqref{eq:virtualw} via propagation through the network with a total communication delay of 200 ms. \icl{Simulations where $\bar{V}_k$ is evaluated via consensus schemes were also carried out} and a similar performance was achieved for small communication delays $< 10$ ms, however the performance deteriorated significantly for larger communication delays.} The switching of the two-level PWM converters causes the DC ripple seen in some of the figures.
\begin{table}[ht]
    \centering
    \caption{Hybrid AC/DC network parameters}
    \begin{tabular}{c|c|c}
     Description & Parameter & Value \\
     \hline
     Bus capacitances & $C_1$, $C_3$, $C_7$ ,$C_9$  & 10 mF \\
     DC load resistances & $R_j, j \in N_{dc}$ & 60 $\Omega$\\
     DC line resistances & $R_{12}$, $R_{23}$, $R_{78}$,$R_{89}$ & 0.01 $\Omega$ \\
     DC line inductances & $L_{12}$, $L_{23}$, $L_{78}$,$L_{89}$ & 0.1 mH \\
     DC switched loads & $P_3$, $P_7$ & 3.6 MW \\
     Rated DC voltage & $V_{dc}^\ast$  & 6000 V\\
     Converter DC capacitances & $C_3$, $C_7$ & 300 mF \\
     Converter AC filter parameters & $R(\Omega), L(mH), C(mF)$ & 0.1, 1, 0.01 \\
     Converter $\omega$ /V ratio & $m$ & 0.002 \\
     AC line resistances & $R_{ij}$ & 0.01 $\Omega$ \\
     AC line inductances & $L_{ij}$ & 0.1 mH \\
     AC load resistances (per phase) & $R_4$, $R_6$ & 60 $\Omega$\\
     AC load active power & $R_5$ & 1 MW\\
     Transformer reactance & $X_5$ & 4\%\\
     Rated AC frequency & $f=\frac{\omega}{2\pi}$ & 50 Hz\\
    \end{tabular}
    \label{tab:acdcparams}
\end{table}
Using the controllers (\ref{eq:droop}), (\ref{eq:ilc}) we show \il{that} the voltages and frequencies \il{converge to equilibrium values} and that the power-sharing is \JDW{close to the optimal values} irrespective of the location of the disturbance.  \JDW{The droop controllers in $\an{\tilde{Q}}$} are set to have gains of ($15.7 kW/\omega$, $10 kW/V$) respectively. In Fig. \ref{fig:b1fv} we show the AC frequency and DC voltage response to step disturbances at time $t = 1s$ and $t = 13s$. The magnitude of the earlier disturbance is 3.6 MW (nominal added demand) located at bus 3 within DC subsystem 1, while the latter disturbance is 3.6 MW reduced demand at bus 7 within DC subsystem 2. The voltages and frequencies of both the ILC buses converge to identical values as expected. Fig. \ref{fig:b1ps} shows that the required power is shared evenly as required.

\begin{figure}[!ht]
    \centering
    \includegraphics[width=0.5\textwidth]{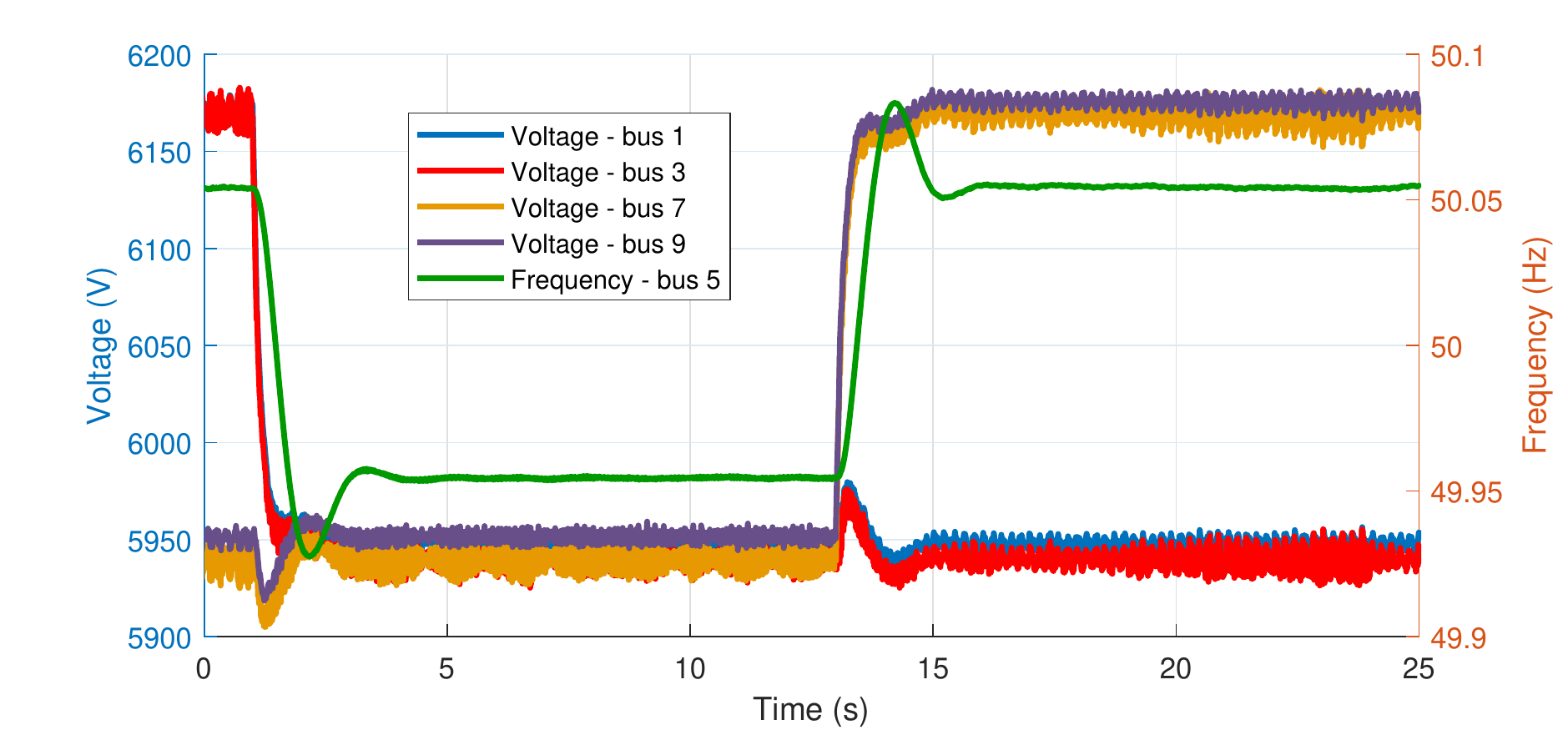}
    \caption{Frequency and voltage response with traditional dual-droop control \eqref{eq:dualdroop}}
    \label{fig:b1fvdd}
\end{figure}

\begin{figure}[!ht]
    \centering
    \includegraphics[width=0.5\textwidth]{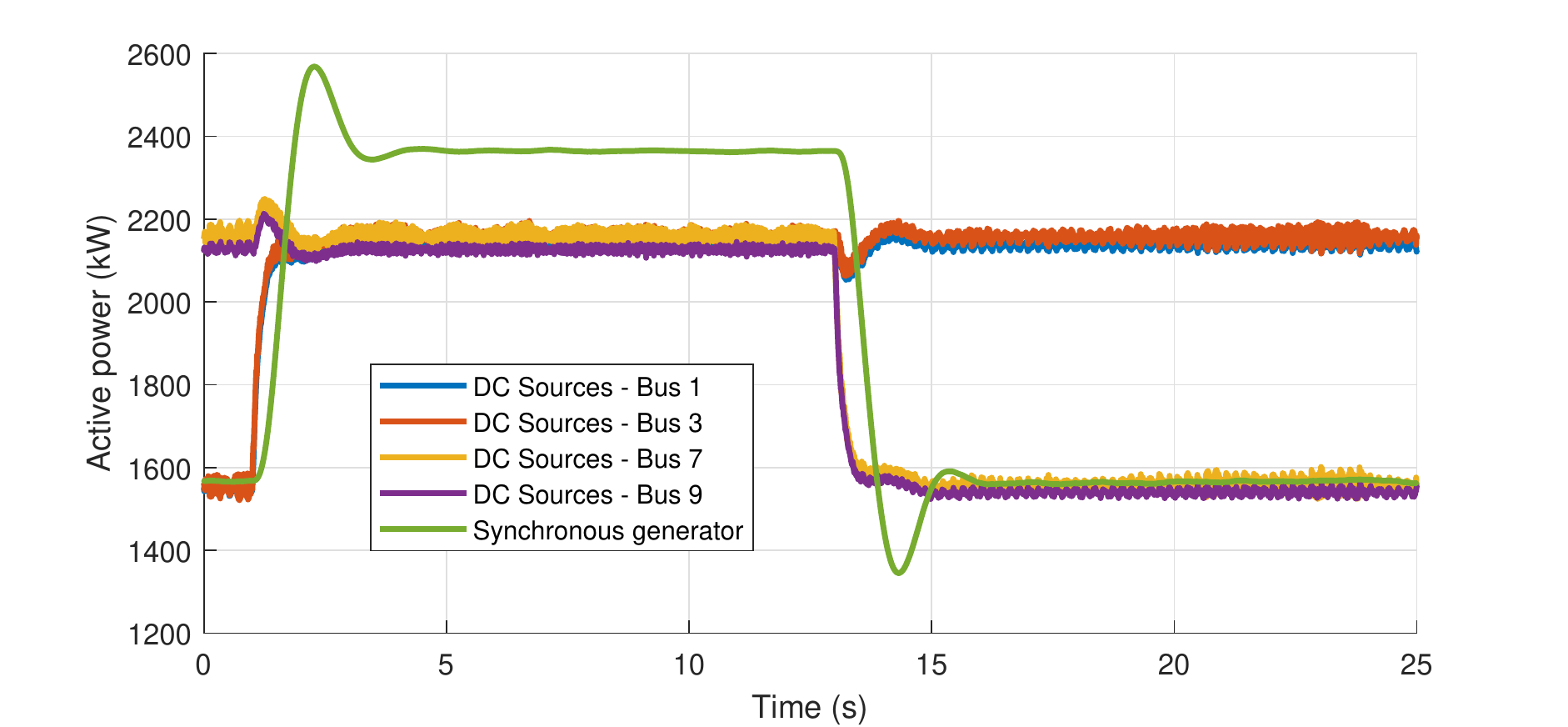}
    \caption{Power-sharing response with traditional dual-droop control \eqref{eq:dualdroop}}
    \label{fig:b1psdd}
\end{figure}

\JDW{For comparison we also study the same hybrid AC/DC network with traditional dual-droop controlled ILCs. The droop gains (4kW / V, 2 MW / (rad/s)) are designed with respect to the respective droop coefficients of the generation in the AC and DC subsystems. Figs. \ref{fig:b1fvdd}-\ref{fig:b1psdd} show that the voltage / frequency deviations at equilibrium as well as the \fv{power sharing} are inferior to the proposed method. Increasing the droop gains decreases the voltage and frequency deviations, however the larger droop gain at the AC and DC sources required to compensate results in worse power-sharing within each DC subsystem, as found in \cite{haileselassie2012}. Furthermore, large disturbances may cause the ILC rating to be exceeded if the droop gains are too large.}
\begin{figure}[ht!]
    \centering
    \includegraphics[width=0.5\textwidth]{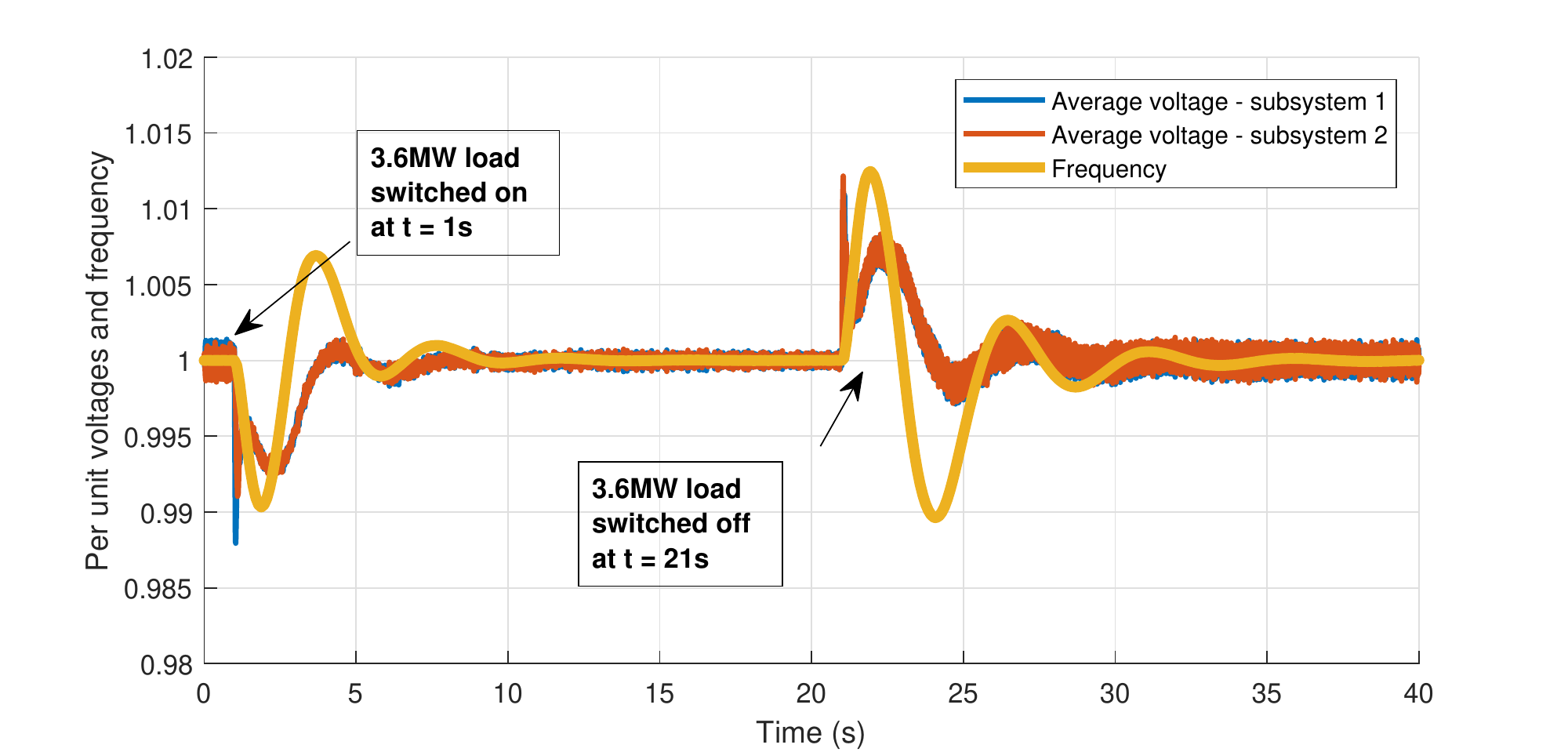}
    \caption{Frequency and voltage response with the \ank{distributed} secondary controller \eqref{eq:secILC},\eqref{eq:sec}}
    \label{fig:b1fv_sec}
\end{figure}
\begin{figure}[ht!]
    \centering
    \includegraphics[width=0.5\textwidth]{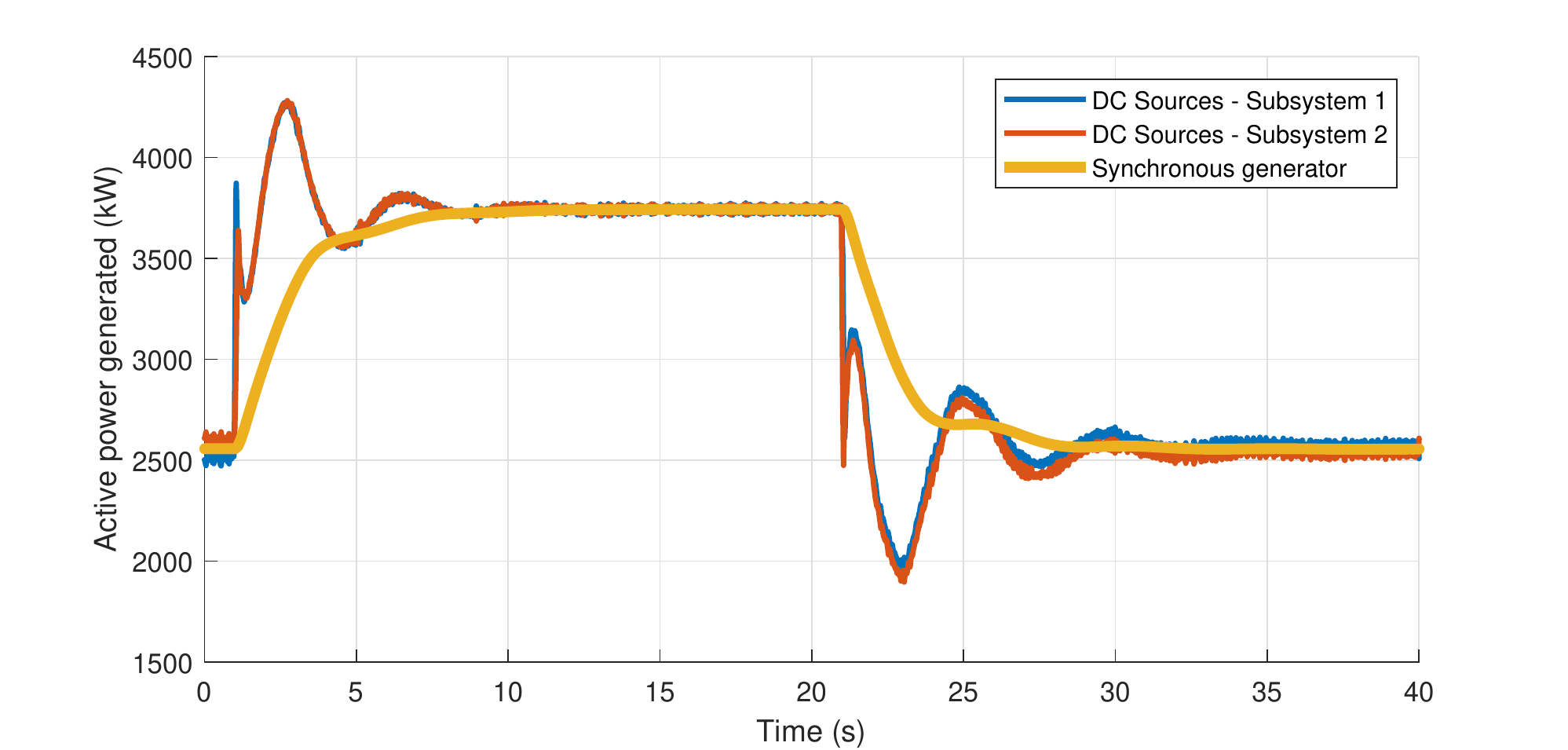}
    \caption{Power-sharing response with the \ank{distributed} secondary controller \eqref{eq:secILC},\eqref{eq:sec}}
    \label{fig:b1ps_sec}
\end{figure}

\JDW{We also consider the performance of the \ak{controller \eqref{eq:secILC},\eqref{eq:sec}} under similar conditions.} Figs. \ref{fig:b1fv_sec}-\ref{fig:b1ps_sec} show that the \ank{distributed} controller regulates the voltages and frequencies of the hybrid AC/DC network to their nominal values while guaranteeing optimal power-sharing at steady-state.

\section{Conclusion} \label{conclusion}
\JDW{We have proposed} a new method \il{for the control of interlinking} converter(s), used in conjunction with traditional droop control to guarantee stability and accurate power sharing in a general hybrid AC/DC network. The stability of the controlled system was proven and it was shown that power-sharing across the AC/DC network is significantly improved compared to dual-droop control.
\il{A secondary control scheme has also been proposed that} guarantees stability while achieving exact optimal power sharing and \il{\ka{that} bus frequencies and the weighted} average of DC voltages \ka{return} to their nominal values \ka{at steady state}. Finally, the proposed algorithms were \il{verified} by simulation and compared to traditional dual-droop control.

\vspace{-1mm}
\section*{Appendix}
\vspace{-1mm}
\il{The appendix includes the proofs of the results presented in the main text. We provide first some notation that will be used within the proofs.}
\JW{Given some column vector $z$ with length $|N|$, we use the subscripted vector $z_{ac}$ to denote \il{the vector that includes the elements of $z$ with indices in $N_{ac}$, i.e. $z_{ac}=[z_{j,j \in N_{ac}}]$}. Likewise the subscripts $z_{Gac}, z_{Xac}, z_{Xdc}, z_{dc}$ denote the vectors \il{that include} those entries for which $j \in N_{ac} \setminus X_{ac}$, $j \in X_{ac}$, $j \in X_{dc}$ and $j \in N_{dc}$ respectively. The following relations therefore hold: $z^T = [z_{ac}^T, z_{dc}^T]$, $z_{ac}^T = [z_{Gac}^T, z_{Xac}^T]$. For convenience we define $\Gamma > 0$ as the diagonal matrix of all $B_{ij}, (i,j) \in E_{ac}$, and $G$ is defined as the conductance matrix of the DC subsystem(s); i.e. the Laplacian weighted by the conductances $G_{ij}$. Finally, summation over all the DC subsystems is denoted by the shorthand $\sum_{k, dc}$.}

\emph{\jw{Proof of Theorem} 1: }We prove our claim in Theorem 1 by finding a suitable Lyapunov function for the system  (\ref{ref:sys}), (\ref{eq:droop}), (\ref{eq:ilc}). Consider the following Lyapunov candidate:
\begin{align}\label{eq:W}
    W\ak{(\eta, \omega^G, V)} =& W_G+W_E+W_V \nonumber\\
    =& \frac{1}{2}(\omega^G-\omega^{G\ast})^TM(\omega^G-\omega^{G\ast})\nonumber\\
    +\mathbf{1}^T\Gamma\int_{\eta*}^{\eta}\sin(\phi)-\sin(\eta^\ast)d\phi &+ \frac{1}{2}m(V-V^\ast)^TC(V-V^\ast).
    \end{align}

The time derivative of the first term is \ak{given by}
\begin{align*}
    \dot{W_G} &= (\omega^G-\omega^{G\ast})^T(p_{Gac}^G-p_{Gac}^L-p_{Gac}^X-p_{Gac}^F-D\omega^G)\\
    &= (\omega^G-\omega^{G\ast})^T(p_{Gac}^G-p_{Gac}^L-p_{Gac}^X-p_{Gac}^F-D\omega^G)\\&\ \ \ +(\omega^X-\omega^{X\ast})^T(p_{Xac}^G-p_{Xac}^L-p_{Xac}^X-p_{Xac}^F)
\end{align*}
\ak{noting that the} second expression follows by adding the term for the converter buses which is equal to zero by \eqref{ref:eq:convbus}. Using the equilibrium conditions \eqref{eqconditions}, \ak{noting that $p_{Gac}^X$ is a zero vector, and rearranging results to}:
\begin{multline}\label{lyap2}
    \dot{W_G} = (\omega-\omega^{*})^T(p_{ac}^G-p_{ac}^{G\ast})
    -\jw{(\omega^X-\omega^{X\ast})^T(p_{Xac}^X-p_{Xac}^{X\ast})}\\
    -(\omega-\omega^{*})^T(p_{ac}^F-p_{ac}^{F\ast})
    -(\omega^G-\omega^{G\ast})^TD(\omega^G-\omega^{G\ast}).
\end{multline}

The time derivative of the second term in \il{\eqref{eq:W}} is, again using the equilibrium conditions (\ref{eqconditions}):
\begin{align*}
    \dot{W_E} &= (\Gamma (sin (\eta)-sin(\eta^\ast)))^T A^T(\omega-\omega^\ast)\\
    &= (\omega-\omega^\ast)^T(p_{ac}^F-p_{ac}^{F\ast}),
\end{align*}
thus canceling the power transfer term in (\ref{lyap2}). The time derivative of \ak{$W_V$ is given by}:
\begin{multline}\label{eq:Wvdot}
    \dot{W_V} = m(V-V^\ast)^T(p_{dc}^G-p_{dc}^{G\ast})
    \jw{-m(V-V^\ast)^T(p_{dc}^X-p_{dc}^{X\ast})}\\
    -m(V-V^\ast)^TG(V-V^\ast).
\end{multline}
Since $G$ is the conductance matrix of the DC graph by definition we have $-m(V-V^\ast)^TG(V-V^\ast) \leq 0$ since $m > 0$. \jdw{Furthermore, using (\ref{eq:ilc}) and noting that $p_{X_x^{dc}}^X=-p_{X_x^{ac}}^X$:}
\jw{
\begin{align*}
    m(V-V^\ast)^T(p_{dc}^X-p_{dc}^{X\ast}) &= (\omega^X - \omega^{X\ast})^T(p_{Xdc}^X-p_{Xdc}^{X\ast})\\
&= -(\omega^X - \omega^{X\ast})^T(p_{Xac}^X-p_{Xac}^{X\ast})
\end{align*}
}
\il{Hence the ILC terms in \eqref{lyap2}, \eqref{eq:Wvdot} can be canceled out.}
 \JW{We also note that converter buses $X_{ac}$ have no frequency-dependent generation nor any damping, \ka{and that} the respective entries of the diagonal matrix $\an{\tilde{Q}}$ are zeros, \ka{while all} other entries of $\an{\tilde{Q}}$ are positive. Therefore, we define \JWW{$\tilde{Q}_G > 0$} as the diagonal matrix with dimension $|N| - |X|$ which includes only the non-zero terms in $\an{\tilde{Q}}$.} Putting it all together \an{and substituting \eqref{eq:droop}}:
\begin{multline}\label{lyapfinalt1}
    \dot{W} \leq -
    \begin{bmatrix}
    (\omega^G-\omega^{G\ast})\\
    m(V-V^\ast)
    \end{bmatrix}^T
    \il{\JWW{\tilde{Q}_G}}
    \begin{bmatrix}
    (\omega^G-\omega^{G\ast})\\
    m(V-V^\ast)
    \end{bmatrix} \\ -(\omega^G-\omega^{G\ast})^TD(\omega^G-\omega^{G\ast}) \leq 0.
\end{multline}

\JW{Using \il{Assumption \ref{as:2}}, $W_E$ has a strict local \il{minimum} 
at $\eta = \eta^\ast$. Likewise $W_G$ and $W_V$ \il{have strict global minima} at $\omega^{G\ast}$ and $V^\ast$ respectively. Thus $W$ has a strict minimum at \jw{$Z^\ast = (\eta^\ast, \omega^{G\ast}, V^\ast)$}. Since \il{$\omega^X$} is uniquely determined by $V$, we can then choose a neighbourhood of \jw{$Z^\ast$} in the coordinates $(\eta, \omega^{G}, V)$. (\ref{lyapfinalt1}) further shows that W is a non-increasing function \il{of time}. Hence the \ak{connected} set \il{$T = \{(\eta, \omega^{G}, V): W \leq \epsilon\}$} for some sufficiently \il{small} $\epsilon > 0$ is compact, forward-invariant and contains \jw{$Z^\ast$}.} \JW{\il{We then apply LaSalle's Theorem,} \il{with $W$ as the Lyapunov-like function,} which states that all trajectories of the system starting from within $T$ converge to the largest invariant set within $T$ \il{that satisfies $\dot W =0$}. Since both \JWW{$\tilde{Q}_G$} and $D$ are positive definite matrices, clearly $\dot{W} = 0$ implies $(\omega^G, V) = (\omega^{G\ast}, V^\ast)$ and therefore $\dot{\omega}^G = \dot{V} = 0$. This in turn implies by (\ref{eqconditions}) that the converter AC-side frequencies $\omega^X = \omega^{X\ast}$. \il{Furthermore, from the equilibrium conditions we deduce that the frequency in each AC subsystem synchronizes to a common value, hence the angle differences $\eta$ converge also to some constant value. Therefore,} by LaSalle's Theorem we have convergence to \ill{the set of equilibrium points} as defined by (\ref{eqconditions}). Finally, choosing \il{$S$ such that it is open, includes \jw{$Z^\ast$}, and $S \subset T$} completes the proof.} \qedsymbol

\emph{\jw{Proof of Theorem} \ref{optim_thm}: }\jw{From the equilibrium conditions \eqref{eqconditions}, $[\omega^T \ mV^T]^T$ \jdw{is arbitrarily close to the vector of equilibrium frequencies} $\mathbf{1}{\omega^\ast}$ as the line resistances become arbitrarily small. This follows from the equilibrium conditions \eqref{eqconditions} and \eqref{eq:pf}, where if the conductances $G_{ij}$ are arbitrarily large, the voltage differences $V_i - V_j$ are arbitrarily small. \jdw{Therefore, to find the power allocation when the line resistances are arbitrarily small, we solve the equilibrium conditions:}}
\begin{align*}
  -\an{\tilde{Q}}\mathbf{1}{\omega}^\ast+\jw{p^G_{nom}}-p^L-p^X-p^F-p^u &= 0\\
    -\mathbf{1}^T\an{\tilde{Q}}\mathbf{1}{\omega}^\ast+\jdw{\mathbf{1}^Tp^G_{nom}}-\mathbf{1}^T(p^L+p^u)-\mathbf{1}^Tp^X-\mathbf{1}^Tp^F &= 0
\end{align*}

Clearly $\mathbf{1}^Tp^F = 0$ in a lossless network and $\mathbf{1}^Tp^X = 0$ for lossless converters. \jdw{Note also that the nominal power generation may be expressed as $p^G_{nom} = -\an{\tilde{Q}}\mathbf{1}\zeta$, where $\zeta$ is a constant.} Hence, solving for $\mathbf{1}{\omega}^\ast$ and substituting into \eqref{eq:droop}:
\begin{align*}
\mathbf{1}{\omega}^\ast =& -\frac{\mathbf{1}\mathbf{1}^T}{\mathbf{1}^T\an{\tilde{Q}}\mathbf{1}}(p^L+p^u) - \jw{\mathbf{1}\zeta}\\
    p^G &= -\an{\tilde{Q}}\mathbf{1}{\omega}^\ast- \jw{\an{\tilde{Q}}\mathbf{1}\zeta} = \jw{\frac{\an{\tilde{Q}}\mathbf{1}\mathbf{1}^T}{\mathbf{1}^T\an{\tilde{Q}}\mathbf{1}}}(p^L+p^u)
\end{align*}
yields the solution (\ref{eq:optimsoln}) to the optimization problem \eqref{eq:optim}. \qedsymbol

\emph{\jw{Proof of Theorem} \ref{optm_thm}: }\il{This is analogous to that of Theorem 2.}
\fv{By premultiplying \eqref{eq:firstpg} by $\mathbf{1}^T$} and noting \eqref{eq:secILC} and the synchronization of frequencies at steady state, it follows that at equilibrium $\hat{\omega} = 0$. The latter shows from \eqref{eq:firstpg} at steady state that
$\xi^* = \mathbf{1}\bar{\xi}$ where $\bar{\xi}$ is the identical equilibrium value of the individual values $\xi_j$ at node $j$.
      Then from the equilibrium conditions (\ref{eqconditions}) we have:
\begin{subequations}
\begin{align}
    \an{\tilde{Q}}\xi^\ast-p^L-p^X-p^F &= 0\\
    \mathbf{1}^T\an{\tilde{Q}}\xi^\ast-\mathbf{1}^Tp^L-\mathbf{1}^Tp^X-\mathbf{1}^Tp^F &= 0
\end{align}
\end{subequations}

Clearly $\mathbf{1}^Tp^F = 0$ in a lossless network and $\mathbf{1}^Tp^X = 0$ for lossless converters. Hence:
\begin{subequations}
\begin{align}
    \mathbf{1}^T\an{\tilde{Q}}\mathbf{1}\bar{\xi}-\mathbf{1}^Tp^L = 0\\
    \xi^\ast=\mathbf{1}\bar{\xi} = \frac{\mathbf{1}\mathbf{1}^Tp^L}{\mathbf{1}^T\an{\tilde{Q}}\mathbf{1}} \label{eq:optimproof}
\end{align}
\end{subequations}

Finally, substituting $p^G = \an{\tilde{Q}}\xi^\ast$ \jw{from \eqref{eq:secpg}} into \eqref{eq:optimproof} yields the solution to the optimization problem \eqref{eq:optim}. \qedsymbol

\emph{\jw{Proof of Theorem} \ref{conv_opt_thm}: }\JW{We use $L = (1, 2, ... |L|)$ to represent the set of nodes of the communication graph.}
\an{For convenience we also write the decomposition of $\tilde{Q}$ into its \il{corresponding} AC and DC sources as $\tilde{Q} = diag(\tilde{Q}_{ac},\tilde{Q}_{dc})$ where $\tilde{Q}_{ac}$ and $\tilde{Q}_{dc}$ are diagonal matrices containing the inverse cost coefficients for the AC and DC generators respectively.}
 We \il{consider the following candidate} Lyapunov function:
\begin{subequations}
\begin{align}
    W &= W_G+W_E+W_V+W_\xi \\
    &= \frac{1}{2}(\omega^{G})^TM\omega^G +\mathbf{1}^T\Gamma\int_{\eta*}^{\eta}sin(\phi)-sin(\eta^\ast)d\phi \nonumber \\ &+ \frac{1}{2}m\bar{V}^T\bar{V} + \frac{1}{2}(\xi-\xi^\ast)^T \ak{T_\xi}(\xi-\xi^\ast)
    \end{align}
\end{subequations}
The time derivatives of $W_G$ are, again adding the term from (\ref{ref:eq:convbus}) and noting that $p_j^X = 0$ for all buses $j \in N_{ac} \setminus X_{ac}$:
\begin{equation*}
\begin{aligned}
    \dot{W}_G = \mbox{ }&(\omega^{G})^T(p_{Gac}^G-p_{Gac}^L-p_{Gac}^F-D\omega^G) \\
    = \mbox{ }&(\omega^{G})^T(p_{Gac}^G-p_{Gac}^L-p_{Gac}^F-D\omega^G)\\
    \mbox{ }&+(\omega^{X})^T(p_{Xac}^G-p_{Xac}^L-p_{Xac}^X-p_{Xac}^F)\\
    =  \mbox{ }&\an{(\omega)^T \tilde{Q}_{ac} (\xi_{ac}-\xi_{ac}^\ast)}-(\omega^{X})^T(p_{Xac}^X-p_{Xac}^{X\ast})\\
    &-\omega^T(p_{ac}^F-p_{ac}^{F\ast})-(\omega^{G})^TD\omega^G\\
\end{aligned}
\end{equation*}
The time derivatives of the other functions comprising $W$ are:
\begin{equation*}
    \begin{aligned}
    \dot{W}_E = \mbox{ }&(\Gamma (sin (\eta)-sin(\eta^\ast)))^T A^T\omega
    = \mbox{ }\omega^T(p_{ac}^F-p_{ac}^{F\ast})\\
    \end{aligned}
\end{equation*}
\begin{equation*}
    \begin{aligned}
    \dot{W_V} =  \mbox{ }&m\bar{V}^T\dot{\bar{V}}\\
    = \mbox{ }&\jw{m\sum_{k, dc}\bar{V}_k\sum_{j\in N_k^{dc}}[(p_j^G-p_j^{G\ast})-(p_j^X-p_j^{X\ast})]}\\
    = \mbox{ }& \jw{
    \ak{m\sum_{k, dc}\bar{V}_k\sum_{j\in N_k^{dc}} \tilde{Q}_{jj} (\xi_j - \xi^*_j)}
    -m\bar{V}^T(p_{Xdc}^X-p_{Xdc}^{X\ast})}\\
     \dot{W}_\xi = \mbox{ }&(\xi-\xi^\ast)^T(-\mathcal{L}(\xi-\xi^\ast) -
     \an{\tilde{Q}}(\hat{\omega}))\\
    = \mbox{ }&-(\xi-\xi^\ast)^T\mathcal{L}(\xi-\xi^\ast)-(\xi-\xi^\ast)^T
     \an{\tilde{Q}}\hat{\omega}\\
    = \mbox{ }&-(\xi-\xi^\ast)^T\mathcal{L}(\xi-\xi^\ast) - \omega^T \an{\tilde{Q}_{ac}}   (\xi_{ac}-\xi_{ac}^\ast) \\\mbox{ }&- \ak{m\sum_{k, dc}\bar{V}_k\sum_{j\in N_k^{dc}} \tilde{Q}_{jj} (\xi_j - \xi^*_j)}
\end{aligned}
\end{equation*}
using (\ref{barv1}) and (\ref{barv2}) to simplify \il{$\dot W_V$}. Clearly $-(\xi-\xi^\ast)^T\mathcal{L}(\xi-\xi^\ast) \leq 0$ from the definition of the Laplacian matrix $\mathcal{L}$. \jdw{We also simplify further by cancelling like terms and thus obtain:}
\begin{multline*}
    \dot{W} \leq -(\omega^G)^TD\omega^G-\jw{(\omega^X)^T(p_{\fv{X}ac}^X-p_{\fv{X}ac}^{X\ast})}\\
    \jw{-m\bar{V}^T(p_{Xdc}^X-p_{Xdc}^{X\ast})}
\end{multline*}
\jdw{Using \eqref{eq:secILC} and noting that $p_{X_x^{dc}}^X=-p_{X_x^{ac}}^X$} we have the cancellation of the last two terms: $    (\omega^{X})^T(p_{Xac}^X-p_{Xac}^{X\ast})=-m\bar{V}^T(p_{Xdc}^X-p_{Xdc}^{X\ast})$. Thus we finally have:
\begin{equation}\label{eq:lyapfinal}
    \dot{W} \leq -(\omega^G)^TD\omega^G \leq 0
\end{equation}
where the damping matrix $D$ is positive definite and can be increased by proportional control of the AC sources. \JW{We then apply LaSalle's Theorem. Clearly, $W$ is minimized at $\eta=\eta^\ast$, $\omega^G = 0$, $\bar{V} = 0$ and $\xi=\xi^\ast$. Therefore we consider the set $T$ which includes $(\eta^\ast, 0, 0,\xi^\ast)$ and is defined by \il{$\{(\eta, \omega^G, \bar{V}, \xi): W \leq \epsilon\}$} for some positive constant $\epsilon$. Since \il{$W$ is non-increasing with time}, \ak{$T$} is a compact, \il{positively invariant set for $\epsilon$ sufficiently small}. LaSalle's Theorem states that trajectories beginning in $T$ converge to the largest invariant set within $T$ for which $\dot{W} = 0$. We therefore examine the equality condition of \eqref{eq:lyapfinal}. \il{$\dot{W} =0$} implies that $\omega^G = 0$, which from the system dynamics (\ref{ref:sys}) implies that $\omega^X = m\bar{V} = 0$, $\bar{V} = 0$, and \il{$\eta$ is constant}. Finally, from (\ref{eq:sec}) if $\hat{\omega} = 0$ then $\mathcal{L}\xi = 0$ which from the definition of the Laplacian communication graph implies that all values $\xi_j, j \in L$ converge to some network-wide \ak{constant value $\bar{\xi}$} and thus $\xi^\ast = \mathbf{1}\bar{\xi}$. Hence the largest invariant set \ak{$\Xi$} within $T$ \il{for which $\dot W =0$} \ak{satisfies $(\eta, \omega^G, \bar{V}, \xi) = (\bar{\eta}, 0, 0, \bar{\xi})$ for constant $\bar{\eta}$ and $\bar{\xi}$. Furthermore, $p^X$ trivially converges from \eqref{ref:eq:convbus}.
To show that within $\Xi$, $V$ \WJ{takes} some constant \icl{value $\hat{V}$ consider} \eqref{eq:dcv} and \eqref{eq:pf} and note that variables $p^G, p^L$ and $p^X$ are constant.
 Then defining $\underline{V}_j = V_j - \hat{V}_j$, it follows that the dynamics of $\underline{V}$ within $\Xi$ satisfy $C \dot{\underline{V}} = -\mathcal{L}_{DC} \underline{V}$ where $\mathcal{L}_{DC}$ is the Laplacian of the graph $(N_{dc}, E_{dc})$, \ka{defined in analogy to \eqref{dfn_laplacian}}.
\icl{It is easy to see that the only invariant set of this linear ODE is}
$\underline{V} \in \text{Im}({\mathbf{1}})$, where $\text{Im}({\mathbf{1}})$ denotes the \ka{image of $\mathbf{1}$},
which together with $\bar{V} = 0$ results to $\underline{V} = 0$.
}
%
\il{Therefore} by LaSalle's theorem the trajectories of the system starting within $T$ converge to \ill{the set of equilibrium points}. \ak{This, in conjunction with Theorem \ref{optm_thm} which suggests that equilibria of \eqref{eq:system},\eqref{ref:sys},\eqref{eq:virtualw},\eqref{eq:secILC},\eqref{eq:sec} are solutions to \eqref{eq:optim} completes} the proof.}~\qedsymbol
%

\balance

%








\end{document}